\documentclass[a4paper,reqno]{amsart}
\usepackage{fullpage}
\usepackage[T1]{fontenc}
\usepackage{lmodern}
\usepackage{amsmath,amssymb,mathtools,booktabs,tikz}
\usetikzlibrary{arrows.meta,calc,decorations.markings,positioning}
\usepackage{microtype}
\usepackage{needspace}
\usepackage{tabularx,array}
\usepackage{xcolor}
\usepackage[hypertexnames=false]{hyperref}
\usepackage{doi}
\usepackage{float}
\hypersetup{colorlinks=true,linkcolor=blue,citecolor=blue,urlcolor=blue}

\newtheorem{thm}{Theorem}[section]
\newtheorem{lem}[thm]{Lemma}
\newtheorem{prop}[thm]{Proposition}

\theoremstyle{definition}

\newtheorem{rem}[thm]{Remark}

\numberwithin{equation}{section}

\newcommand{\unorm}[1]{%
	\left|\mkern-1.5mu\left|\mkern-1.5mu\left|#1\right|
	\mkern-1.5mu\right|\mkern-1.5mu\right|%
}
\newcommand{\M}{\mathbb M}
\newcommand{\C}{\mathbb C}
\newcommand{\R}{\mathbb R}
\newcommand{\N}{\mathbb N}
\newcommand{\Tr}{\operatorname{Tr}}
\newcommand{\diag}{\operatorname{diag}}

\newcommand{\cD}{\mathcal{D}}

\begin{document}
	
	\title[Clarkson--McCarthy type inequalities, part II]
	{Clarkson--McCarthy type inequalities, part II: refinements, unitary orbits, unitarily invariant norms and $\ell_s$--$\ell_s$ estimates}

	\author[T.~Zhang]{Teng Zhang}
	
	\address{School of Mathematics and Statistics, Xi'an Jiaotong University, Xi'an 710049, P. R. China}
	\email{teng.zhang@stu.xjtu.edu.cn}
	

	
	\begin{abstract}
In previous work \cite{Zha26a}, the author established Clarkson--McCarthy type
inequalities for several classes of operator families, including
$\ell_p$--$\ell_p$ estimates and mixed $\ell_q$--$\ell_p$ estimates, by means of
operator square-sum identities. The present paper has two main aims. First,
using the tracial joint spectral measures recently introduced by
O.~Hein\"avaara \cite{Hei25}, we obtain several refinements of the
McCarthy--Clarkson inequalities, one of which confirms a conjecture of
Audenaert \cite{AK17}. Second, combining tools from matrix analysis with
complex interpolation theory, we establish Clarkson--McCarthy type inequalities
for unitary orbits, unitarily invariant norms, and $\ell_s$--$\ell_s$ estimates,
again in the framework of operator square-sum identities.
	\end{abstract}

	\maketitle
	\section{Introduction} 
	
	For a compact operator $T$ on a fixed Hilbert space and $0<p<\infty$,
	we write $\|T\|_p=(\Tr|T|^p)^{1/p}$, where $|T|=(T^*T)^{1/2}$.
	The Schatten class $S_p$ consists of the compact operators on this
	space for which $\|T\|_p<\infty$. On $S_p$, this functional is a norm
	for $p\ge1$ and a quasi-norm for $0<p<1$.
	We denote the operator norm by $\|T\|_\infty$.
	
	In 1967, McCarthy \cite{McC67} obtained two uniform convexity
	inequalities for Schatten trace ideals $S_p$, now known as the
	\emph{Clarkson--McCarthy inequalities}.  The two inequalities needed
	below can be stated as follows.  Let $A,B\in S_p$.  Then
	\begin{align}
		\|A+B\|_p^p+\|A-B\|_p^p
		&\le 2^{p-1}\bigl(\|A\|_p^p+\|B\|_p^p\bigr),
		\qquad 2\le p<\infty, \label{eq:McC1}\\
		\|A+B\|_p^q+\|A-B\|_p^q
		&\le 2\bigl(\|A\|_p^p+\|B\|_p^p\bigr)^{q/p},
		\qquad 1<p\le 2, \label{eq:McC-mixed}
	\end{align}
	where $q$ denotes the conjugate exponent of $p$, i.e.\
	$1/p+1/q=1$.  Moreover, \eqref{eq:McC1} holds with the reverse
	inequality for $0<p\le2$, and \eqref{eq:McC-mixed} holds with the
	reverse inequality for $2\le p<\infty$.
	These two inequalities show that $S_p$ is $q$-uniformly convex for $1 < p \le 2$ and $p$-uniformly convex for $2\le p<\infty$ . However, McCarthy was not the first to prove \eqref{eq:McC1}; this inequality already appeared in work of Dixmier \cite{Dix53}. McCarthy's original proof of \eqref{eq:McC-mixed} contains a gap; see \cite[p.~299]{FK86}. A proof by complex interpolation is given in \cite[Theorem~5.3]{FK86}.  In the commutative setting, Clarkson's inequalities~\cite{Cla36} in $L_p$ were sharpened by Hanner \cite{Han56}; Hanner's inequality can be  used to compute the optimal moduli of convexity and smoothness of $L_p$ and $\ell_p$.   Ball, Carlen, and Lieb \cite{BCL94} then proved the sharp uniform convexity and smoothness estimates for $S_p$, and in the same work formulated the following noncommutative Hanner's inequality problem for Schatten classes: for $A,B\in S_p$, 
	\begin{equation}\label{eq:Hanner}
		\|A+B\|_p^p+\|A-B\|_p^p
		\le
		\bigl(\|A\|_p+\|B\|_p\bigr)^p
		+
		\bigl|\|A\|_p-\|B\|_p\bigr|^p, \quad p\ge 2.
	\end{equation} 
	The inequality is reversed for $1\le p\le 2$. For earlier work on \eqref{eq:Hanner}, see \cite{BCL94,Kin03,CL06,Cha21a,Hei24}. Hein\"avaara~\cite{Hei25} proved Hanner's inequality in $S_p$ throughout $1\le p<\infty$ using \emph{tracial joint spectral measures}. More precisely, Hein\"{a}vaara \cite{Hei25} proved that for fixed $A,B\in S_p$, the real linear span of $A$ and $B$ is isometric to a subspace of $L_p(\mu)$ for some positive measure $\mu$.

	Related extensions concern several operators
	\cite{BK04,HK08,Kec19,Zha25,Zha26a}, unitary orbits \cite{BL20},
	unitarily invariant norms \cite{BH88}, and the quadratic symmetric
	modulus \cite{Zha26b}.

	\subsection{Refinements of Clarkson--McCarthy inequality \texorpdfstring{\eqref{eq:McC1}}{(1.1)}}
	We first sharpen \eqref{eq:McC1} for all $p\ge1$.
	The following result improves the coefficient in
	\cite[Theorem~3.2]{AHK18} and removes the additional positivity
	assumptions imposed there in the intermediate exponent ranges.
	\begin{thm}\label{thm:sharp-clarkson-matrix-upper}
		Let $A,B\in S_p$. Then for $p\ge 2$,
		\begin{equation*}
			\|A+B\|_p^p+\|A-B\|_p^p
			\le
			2^{p-1}\bigl(\|A\|_p^p+\|B\|_p^p\bigr)
			-
			\bigl(2^{p-1}-2\bigr)
			\bigl|\|A\|_p-\|B\|_p\bigr|^p .
		\end{equation*}
		For $1\le p\le2$, the inequality is reversed.
		Moreover, the coefficient pair
		$
		\bigl(2^{p-1},\,2^{p-1}-2\bigr)
		$
		is optimal: for $p\ge2$, the first coefficient cannot be decreased
		and, with it fixed, the second cannot be increased. For $1\le p\le2$,
		the first coefficient cannot be increased and, with it fixed, the
		second cannot be decreased.
	\end{thm}
	Let $\mathbb{M}_n$ be the algebra of all $n\times n$ complex matrices, $\mathbb{H}_n$ the real vector space of all Hermitian matrices, $\mathbb{P}_n$ the cone of positive semidefinite matrices, and $\mathbb{P}_n^+$ the cone of positive definite matrices. For $X,Y\in \mathbb{H}_n$, we write $X\le Y$ if $Y-X\in\mathbb{P}_n$.
	Let $A,B\in \mathbb P_n^+$. Applying \eqref{eq:McC1} with exponent $2p$ to $
	\frac12
	\begin{pmatrix}
		A^{1/2}&B^{1/2}\\
		0&0
	\end{pmatrix},
	\frac12
	\begin{pmatrix}
		A^{1/2}&-B^{1/2}\\
		0&0
	\end{pmatrix}, 
	$
	we obtain the positive semidefinite trace form
	\begin{equation}\label{eq:trace-superadd}
		\Tr(A+B)^p
		\ge
		\Tr A^p+\Tr B^p,
		\qquad  p\ge1.
	\end{equation}
	For $0<p\le1$, the same argument,
	using the reverse form of \eqref{eq:McC1}, gives the reverse inequality
	$
	\Tr(A+B)^p
	\le
	\Tr A^p+\Tr B^p .
	$
	In fact, for $p>0$ both directions are also direct consequences of
	Rotfel'd's trace inequality \cite{Rot69}.  For negative powers, the
	reverse estimate follows instead from Weyl's eigenvalue monotonicity.  If
	$A,B\in\mathbb P_n^+$ and $p<0$, then $A+B\ge A$, and hence
	$
	\lambda_j(A+B)\ge \lambda_j(A), j=1,\ldots,n,
	$ where $\lambda_j(\cdot)$ is the $j$th largest eigenvalue.
	Since $t^p$ is decreasing on $(0,\infty)$ for $p<0$, we obtain
	$
	\Tr(A+B)^p
	\le
	\Tr A^p
	\le
	\Tr A^p+\Tr B^p .
	$
	Thus the reverse inequality of \eqref{eq:trace-superadd} holds for all $p\le1$.

	We also give a proof of the power-trace inequality conjecture proposed by
	Audenaert in \cite[Section~3]{AK17}, covering all real exponents.
	Audenaert proved the cases $p\le-2$ and $0<p\le3$ by an integral
	representation of power functions~\cite{Aud12}.
	Related block-matrix inequalities were stated in \cite{Aud06}.
Our proof treats positive exponents
	using tracial joint spectral measures and negative exponents using
	geometric means and log-majorization.
	\begin{thm}\label{thm:main}
		Let $A,B\in\mathbb P_n^+$.  For $0\le p\le1$ and for $p\ge2$, we have
		\begin{equation}\label{eq:McC_refine}
			\Tr(A+B)^p
			\ge 
			\Tr A^p+\Tr B^p
			+
			(2^p-2)\Tr(A^{1/2}BA^{1/2})^{p/2}.
		\end{equation}
		For $p\le 0$ and for $1\le p\le 2$, the inequality is reversed.
		At $p=0,1,2$ equality holds for all $A,B$.  For every $p$ in the
		above ranges, equality is also attained when $A=B$. Moreover, the coefficient triple
		$
		\bigl(1,1,2^p-2\bigr)
		$
		is optimal in the following componentwise sense.  In the lower-bound ranges
		$0\le p\le1$ and $p\ge2$, none of the three coefficients can be
		increased while the other two are kept fixed. In the upper-bound
		ranges $p\le 0$ and $1\le p\le2$, none of the three coefficients can be
		decreased while the other two are kept fixed.  
	\end{thm}
	
	Theorem~\ref{thm:main} has different interpretations in different
	ranges of $p$.  For $p\ge2$, since $2^p-2\ge0$ and the inequality \eqref{eq:McC_refine} has the
	same direction as \eqref{eq:trace-superadd}, it gives a genuine
	refinement of \eqref{eq:trace-superadd}.  For $1\le p\le2$, it gives a
	complementary upper bound to \eqref{eq:trace-superadd}.  For
	$0<p\le1$, it gives a complementary lower bound to the reversed
	Rotfel'd inequality.  For $p<0$, it strengthens the elementary reverse
	estimate
	$
	\Tr(A+B)^p\le \Tr A^p+\Tr B^p .
	$

	The proof of Theorem~\ref{thm:main} separates naturally into two regimes.  Positive exponents are handled by a tracial joint spectral measure, where the problem becomes a scalar inequality on the cone $a\ge |b|$.  Negative exponents are handled by geometric means and log-majorization.  
	
	\begin{rem}
		We have learned that Li and Zhou~\cite[Theorem~1.1(i)]{LZ26}
		also proved the inequality in Theorem~\ref{thm:main}
		for all negative exponents. 
	\end{rem}
	\subsection{Hanner's inequality and norm compression below one}
	\label{subsec:below-one}
	Although Hanner's inequality~\eqref{eq:Hanner} holds throughout the Banach range,
	no comparison with a fixed direction survives for arbitrary
	matrices when $0<p<1$.
	
	\begin{prop}\label{prop:hanner-no-direction}
		For every $0<p<1$, each of the two strict inequalities between
		\[
		\|A+B\|_p^p+\|A-B\|_p^p
		\quad\text{and}\quad
		(\|A\|_p+\|B\|_p)^p+|\|A\|_p-\|B\|_p|^p
		\]
		occurs for a pair of commuting positive semidefinite diagonal matrices.
	\end{prop}
	
	Positivity of both $A+B$ and $A-B$ changes the conclusion.
	The following reverse Hanner inequality was stated by
	Chayes~\cite[Theorem~1.2]{Cha21b}. However, her proof contains a gap; specifically, the coefficient comparison in
	\cite[(4.20)]{Cha21b} has the wrong direction: for $0<p<1$,
	the factor $p(p-1)\cdots(p-2k+1)$ is negative for every $k\ge1$.
	We give a proof from the block-matrix inequality below, independent
	of that argument.
	
	\begin{thm}\label{thm:hanner-positive}
		Let $0<p<1$ and let $A,B\in\mathbb H_n$ satisfy $A+B\ge0$ and
		$A-B\ge0$. Then
		\begin{equation}\label{eq:hanner-positive}
			\|A+B\|_p^p+\|A-B\|_p^p
			\le(\|A\|_p+\|B\|_p)^p+|\|A\|_p-\|B\|_p|^p.
		\end{equation}
	\end{thm}
	
	The connection with positive block matrices goes back to
	King~\cite[Theorem~1]{Kin03}. For
	$M=\bigl(\begin{smallmatrix}X&Y\\Y^*&Z\end{smallmatrix}\bigr)\ge0$,
	write
	\[
	m_p=\begin{pmatrix}\|X\|_p&\|Y\|_p\\
		\|Y\|_p&\|Z\|_p\end{pmatrix}.
	\]
	King proved $\|M\|_p\ge\|m_p\|_p$ for $1\le p\le2$ and the
	reverse inequality for $p\ge2$. The special choice $X=Z$ and
	$Y=Y^*$ gives the positive-matrix case of Hanner's inequality (Theorem~\ref{thm:hanner-positive}).
	Our next result establishes the reverse norm compression inequality
	throughout $0<p<1$.
	
	\begin{thm}\label{thm:reverse-king-half-one}
		Let $X,Z\in\mathbb P_n$, $Y\in\M_n$, and suppose that
		$M=\bigl(\begin{smallmatrix}X&Y\\Y^*&Z\end{smallmatrix}\bigr)\ge0$.
		For every $0<p<1$, its scalar compression $m_p$ is positive
		semidefinite and
		\begin{equation}\label{eq:reverse-king}
			\left\|\begin{pmatrix}X&Y\\Y^*&Z\end{pmatrix}\right\|_p
			\le
			\left\|\begin{pmatrix}\|X\|_p&\|Y\|_p\\
				\|Y\|_p&\|Z\|_p\end{pmatrix}\right\|_p.
		\end{equation}
	\end{thm}
	
	The proof starts with the reverse H\"older inequality and its
	negative-exponent variational formula. A geometric-mean construction
	proves \eqref{eq:reverse-king} for $1/2\le p<1$.
	A scalar concavity estimate then transfers the result from
	$p/(1-p)$ to $p$, and induction covers all positive exponents
	below one. 
	
	\subsection{Several operators and mixed norms}
	\label{subsec:families}
	
	The remaining results concern matrix families with equal square
	sums and operator families related by a scalar isometry. They extend
	the square-sum approach in
	Zhang~\cite{Zha26a}, in the setting of the several-operator
	inequalities of Bhatia and Kittaneh~\cite{BK04} and Hirzallah and
	Kittaneh~\cite{HK08}. Unitary orbit versions of the two-operator
	Clarkson--McCarthy inequalities were studied by Bourin and
	Lee~\cite{BL20}; unitarily invariant norm versions were obtained
	by Bhatia and Holbrook~\cite{BH88}.
	
	The next two theorems require only equality of the square sums.
	This condition holds, in particular, when
	$B_i=\sum_{j=1}^n u_{ij}A_j$ and
	$U=(u_{ij})\in\M_{m,n}(\C)$ satisfies $U^*U=I_n$;
	see \eqref{eq:isometry-square-sum}.
	Here $\M_{m,n}(\C)$ denotes the space of complex $m\times n$
	matrices, and such a matrix $U$ is called an isometry.
	
	\begin{thm}\label{thm:unitary-orbit}
		Let $A_1,\ldots,A_n,B_1,\ldots,B_m\in\M_d$ satisfy
		\[
		\sum_{i=1}^m|B_i|^2=\sum_{j=1}^n|A_j|^2.
		\]
		For $p\ge2$, there are unitaries $U_i,V_j\in\M_d$ such that
		\begin{align}
			\sum_{i=1}^mU_i|B_i|^pU_i^*
			&\le n^{p/2-1}\sum_{j=1}^n|A_j|^p,
			\label{eq:orbit-forward}\\
			\sum_{j=1}^nV_j|A_j|^pV_j^*
			&\le m^{p/2-1}\sum_{i=1}^m|B_i|^p.
			\label{eq:orbit-backward}
		\end{align}
		For $0<p\le2$, there are unitaries for which both inequalities reverse.
	\end{thm}
	
	Theorem~\ref{thm:unitary-orbit} proves the author's unitary orbit
	conjecture in \cite{Zha25}. Indeed, the identity
	\[
	\left|\sum_{j=1}^n A_j\right|^2+
	\sum_{1\le j<k\le n}|A_j-A_k|^2
	=n\sum_{j=1}^n|A_j|^2
	\]
	allows us to apply \eqref{eq:orbit-forward} to the family
	$(\sqrt n A_j)_{j=1}^n$ and the family consisting of
	$\sum_{j=1}^n A_j$ and all the differences $A_j-A_k$, $j<k$.
	Thus, for $p\ge2$, there are unitaries $W_0,W_{jk}\in\M_d$ such that
	\[
	W_0\left|\sum_{j=1}^n A_j\right|^pW_0^*
	+\sum_{1\le j<k\le n}W_{jk}|A_j-A_k|^pW_{jk}^*
	\le n^{p-1}\sum_{j=1}^n|A_j|^p.
	\]
	For $0<p\le2$, suitable unitaries give the reverse inequality.
	Theorem~\ref{thm:uin} gives the corresponding norm estimates.
	
	We write $\unorm{\cdot}$ for an arbitrary unitarily invariant norm
	on $\M_d$, so that $\unorm{VAT}=\unorm{A}$ for all unitary
	$V,T\in\M_d$.
	
	\Needspace{10\baselineskip}
	\begin{thm}\label{thm:uin}
		Let $A_1,\ldots,A_n,B_1,\ldots,B_m\in\M_d$ satisfy
	\[
		\sum_{i=1}^m|B_i|^2=\sum_{j=1}^n|A_j|^2.
		\] For $p\ge2$,
		\begin{equation}\label{eq:uin-family}
			m^{1-p/2}\unorm{\sum_{j=1}^n|A_j|^p}
			\le\unorm{\sum_{i=1}^m|B_i|^p}
			\le n^{p/2-1}\unorm{\sum_{j=1}^n|A_j|^p}.
		\end{equation}
		Both inequalities reverse for $0<p\le2$.
	\end{thm}
	
	Taking traces in \eqref{eq:orbit-forward} and
	\eqref{eq:orbit-backward}, or taking the Schatten $1$-norm in
	\eqref{eq:uin-family}, gives
	\[
	m^{1-p/2}\sum_{j=1}^n\|A_j\|_p^p
	\le\sum_{i=1}^m\|B_i\|_p^p
	\le n^{p/2-1}\sum_{j=1}^n\|A_j\|_p^p,
	\qquad p\ge2,
	\]
	with both inequalities reversed for $0<p\le2$.
	These are the finite-dimensional inequalities of
	\cite[Theorem~1.8]{Zha26a}; their extension to Schatten classes
	is explained in Section~\ref{sec:families}.

	For families related by a scalar isometry, keeping the Schatten
	exponent fixed and interpolating in the outer sequence exponent
	gives the following $\ell_s(S_p)$ estimates.
	For a finite tuple, we adopt the notation
	\[
	\|(A_j)_{j=1}^n\|_{\ell_s^n(S_p)}=
	\begin{cases}
		(\sum_{j=1}^n\|A_j\|_p^s)^{1/s},&1\le s<\infty,\\
		\max_{1\le j\le n}\|A_j\|_p,&s=\infty.
	\end{cases}
	\]
	With $1/\infty=0$, define
	\begin{equation}\label{eq:mixed-constants}
		C_{p,s}(m,n)=
		\begin{cases}
			m^{1/s-1/2},&1\le s\le p,\quad1\le p\le2,\\[1mm]
			m^{(1-p/2)/s}n^{(1-p/s)/2},
			&p\le s\le\infty,\quad1\le p\le2,\\[1mm]
			m^{\frac{p/s-1}{2(p-1)}}n^{\frac{(p-2)(1-1/s)}{2(p-1)}},
			&1\le s\le p,\quad2\le p<\infty,\\[1mm]
			n^{1/2-1/s},&p\le s\le\infty,\quad2\le p<\infty.
		\end{cases}
	\end{equation}
	The formulas agree at $p=2$ and $s=p$.
	
	\begin{thm}\label{thm:schatten_p_q_norms}
		Let $1\le p<\infty$, $1\le s\le\infty$, and $A_1,\ldots,A_n\in S_p$.
		Let $U\in\M_{m,n}(\C)$ satisfy $U^*U=I_n$, and put
		$B_i=\sum_{j=1}^n u_{ij}A_j$. Then
		\begin{align}
			\|(B_i)_{i=1}^m\|_{\ell_s^m(S_p)}
			&\le C_{p,s}(m,n)\|(A_j)_{j=1}^n\|_{\ell_s^n(S_p)},
			\label{eq:ell-s-forward}\\
			\|(A_j)_{j=1}^n\|_{\ell_s^n(S_p)}
			&\le C_{p,s}(n,m)\|(B_i)_{i=1}^m\|_{\ell_s^m(S_p)}.
			\label{eq:ell-s-backward}
		\end{align}
	\end{thm}
	
	In fact, \eqref{eq:ell-s-forward} holds whenever $U$ is a scalar
	contraction, $U^*U\le I_n$. Establishing this stronger forward
	statement makes it possible to apply the same estimate to $U^*$
	and prove \eqref{eq:ell-s-backward} in the isometric case.
	
	For scalar isometries and $1\le p<\infty$,
	Theorem~\ref{thm:schatten_p_q_norms} extends the
	$\ell_p(S_p)$ estimates in \cite[Theorem~1.8]{Zha26a}
	to $\ell_s(S_p)$ for every $1\le s\le\infty$.
	At $s=p$, the constants agree with those in that result.
	In certain ranges of $p$ and $s$, interpolation also improves
	the bounds obtained from that result by elementary comparisons
	of finite sequence norms.

	\medskip
	\noindent\textbf{Organization of this paper.}
		The paper is organized as follows. Sections~\ref{sec:sharp-clarkson}
	and~\ref{sec:trace-refinement} prove the two refinements in
	Section~1.1. Section~\ref{sec:hanner-below-one} proves the assertions
	about Hanner's inequality below one, using
	Theorem~\ref{thm:reverse-king-half-one} for the positive case.
	Section~\ref{sec:reverse-king} proves that norm compression theorem.
	Section~\ref{sec:families} treats unitarily invariant norms, unitary
	orbits, and mixed norms in that order.

	\medskip
	\noindent\textbf{Acknowledgments and AI tools disclosure. }
	Teng Zhang is supported by the China Scholarship Council, the Young Elite Scientists Sponsorship Program for PhD Students of the China Association for Science and Technology, and the Fundamental Research Funds for the Central Universities at Xi'an Jiaotong University (Grant No.~xzy022024045).
	
	AI tools were used during the preparation of this paper to aid with proofreading, editing and content critique. All mathematical ideas are due to the author.
	
	\section{Proof of Theorem~\ref{thm:sharp-clarkson-matrix-upper}}
	\label{sec:sharp-clarkson}
	In this section, we give a proof of Theorem~\ref{thm:sharp-clarkson-matrix-upper}.
	\begin{proof}[Proof of Theorem~\ref{thm:sharp-clarkson-matrix-upper}]
		Put $a=\|A\|_p$ and $b=\|B\|_p$, and assume, without loss of generality, that $a\ge b$.  By 
		Hanner's inequality  \eqref{eq:Hanner}, for $p\ge2$, we have
		\[
		\|A+B\|_p^p+\|A-B\|_p^p
		\le
		(a+b)^p+(a-b)^p .
		\]
		Now, we claim that
		\begin{equation}\label{eq:sharp-clarkson-scalar}
			(a+b)^p+(a-b)^p
			\le
			2^{p-1}(a^p+b^p)
			-
			(2^{p-1}-2)(a-b)^p .
		\end{equation}
		Indeed, set $x=a+b$ and $y=a-b$.  Then $x\ge y\ge0$, and
		\eqref{eq:sharp-clarkson-scalar} is equivalent to
		\[
		(x+y)^p+(x-y)^p-2x^p-(2^p-2)y^p\ge0 .
		\]
		If $x=0$ there is nothing to prove.  Otherwise put $t=y/x\in[0,1]$ and
		\[
		\phi_p(t)=(1+t)^p+(1-t)^p-2-(2^p-2)t^p .
		\]
		Since $\phi_p(0)=\phi_p(1)=0$, only $0<t<1$ needs consideration.
		For $0<t\le1$, define
		\[
		g_p(t)=\frac{(1+t)^p+(1-t)^p-2}{t^p}.
		\]
		For $0<t<1$, differentiation gives
		\[
		g_p'(t)
		=
		\frac{p\bigl(2-(1+t)^{p-1}-(1-t)^{p-1}\bigr)}{t^{p+1}} .
		\]
		For $p\ge2$, the function $u\mapsto u^{p-1}$ is convex on
		$(0,\infty)$; hence
		\[
		(1+t)^{p-1}+(1-t)^{p-1}\ge2,
		\]
		so that $g_p'(t)\le0$.  Therefore
		\[
		g_p(t)\ge g_p(1)=2^p-2,
		\]
		which proves \eqref{eq:sharp-clarkson-scalar}.  Combining this scalar
		estimate with Hanner's inequality proves the desired estimate for
		$p\ge2$.
		
		For $1\le p\le2$,  Hanner's inequality  \eqref{eq:Hanner} is reversed, and the
		same scalar argument is reversed because $u\mapsto u^{p-1}$ is concave
		on $(0,\infty)$ for $1<p<2$.  Hence the asserted reverse inequality
		follows for $1<p<2$. At $p=1$, the asserted inequality is
		$\|A+B\|_1+\|A-B\|_1\ge2\max\{\|A\|_1,\|B\|_1\}$,
		which follows by applying the triangle inequality to $2A$ and $2B$.
		At $p=2$, equality is the parallelogram identity.
		
		It remains to check sharpness. First let $p\ge2$. If the coefficient of
		$\|A\|_p^p+\|B\|_p^p$ were replaced by a smaller constant $C$, then
		taking $A=B\ne0$ would give
		$
		2^p\|A\|_p^p\le 2C\|A\|_p^p,
		$
		so $C\ge2^{p-1}$.  Thus $2^{p-1}$ is sharp.  With
		$C=2^{p-1}$ fixed, taking $B=0$ gives
		$
		2\|A\|_p^p
		\le
		\bigl(2^{p-1}-D\bigr)\|A\|_p^p,
		$
		so $D\le2^{p-1}-2$.  Hence the coefficient $2^{p-1}-2$ is also sharp.
		For $1\le p\le2$, the same substitutions in the reversed inequality
		give $C\le2^{p-1}$ and, with $C=2^{p-1}$, $D\ge2^{p-1}-2$.
		This proves the optimality assertions in both ranges.
	\end{proof}

	\section{Proof of Theorem~\ref{thm:main}}
	\label{sec:trace-refinement}
	In this section, we give a proof of Theorem~\ref{thm:main}.
	\subsection{Tracial joint spectral measure}
	We use the \emph{tracial joint spectral measure} introduced by
	Hein\"avaara~\cite{Hei25} in his proof of Hanner's inequality
	\eqref{eq:Hanner}.
	\begin{lem}[{\cite[Theorem~1.3]{Hei25}}]\label{thm:heinavaara}
		Let $A,B\in \mathbb{H}_n$.  There exists a positive measure $\mu_{A,B}$ on $\R^2$, called \emph{tracial joint spectral measure of $A,B$}.  Fix any real-valued measurable function $f$ on $\R$ such that $f(0)=0$, and for any $M>0$,
		\[
		\int_{-M}^{M}\left|\frac{f(t)}{t}\right|\,dt<\infty.
		\]
		Define a function $H(f):\R\to\R$ by
		\[
		H(f)(s)=\int_0^1\frac{1-t}{t}f(st)\,dt.
		\]
		Then for all $x,y\in\R$, we have
		\begin{equation}\label{eq:tjsm-general}
			\Tr H(f)(xA+yB)=\int_{\R^2} f(ax+by)\,d\mu_{A,B}(a,b).
		\end{equation}
		Consequently, for every $p>0$,
		\begin{equation}\label{eq:tjsm-p}
			\Tr |xA+yB|^p=p(p+1)\int_{\R^2}|ax+by|^p\,d\mu_{A,B}(a,b).
		\end{equation}
	\end{lem}
	
	The following support lemma  is an elementary consequence of Lemma~\ref{thm:heinavaara}. Here and below, a property is said to hold $\mu$-almost everywhere
	($\mu$-a.e.) if the set on which it fails has $\mu$-measure zero.
	
	\begin{lem}\label{lem:support}
		Let $D,O\in \mathbb{H}_n$ satisfy
		$
		D\pm O\ge0.
		$
		Let $\mu=\mu_{D,O}$ be the tracial joint spectral measure of $D,O$ determined in Lemma~\ref{thm:heinavaara}.  Then
		\begin{equation*}
			a\ge |b|\qquad \mu\hbox{-a.e.}
		\end{equation*}
	\end{lem}
	
	\begin{proof}
		Write $s_+=\max\{s,0\}$. Fix $p>0$ and take
		$
		f(s)=(-s)_+^p
		$ in Lemma~\ref{thm:heinavaara}.
		Then
		\[
		H(f)(s)=\int_0^1\frac{1-t}{t}(-st)_+^p\,dt
		=\frac{(-s)_+^p}{p(p+1)}.
		\]
		Since $D+O\ge0$, we have $H(f)(D+O)=0$.  Applying \eqref{eq:tjsm-general} with $(x,y)=(1,1)$ gives
		\[
		0=\int_{\R^2}(-a-b)_+^p\,d\mu(a,b).
		\]
		The measure $\mu$ is positive, hence $a+b\ge0$ $\mu$-almost everywhere.  Applying the same argument to $D-O\ge0$, namely with $(x,y)=(1,-1)$, gives $a-b\ge0$ $\mu$-almost everywhere.  Thus $a\ge |b|$ $\mu$-almost everywhere.
	\end{proof}
	
	The entire sign pattern for positive exponents is encoded in the following scalar inequality.
	
	\begin{lem}\label{lem:scalar}
		Let $a,b\in\R$ satisfy $a\ge |b|$.  For $p>0$, define
		\begin{equation}\label{eq:scalar-phi}
			\Phi_p(a,b):=(a+b)^p+(a-b)^p-2a^p-(2^p-2)|b|^p.
		\end{equation}
		Then
		\[
		\Phi_p(a,b)
		\begin{cases}
			\ge0, & 0<p<1\ \hbox{or}\ p>2,\smallskip\\
			=0, & p=1,2,\smallskip\\
			\le0, & 1<p<2.
		\end{cases}
		\]
	\end{lem}
	
	\begin{proof}
		If $a=0$, then $b=0$ and there is nothing to prove.  Assume $a>0$ and write
		$
		t=|b|/a\in[0,1].
		$
		It is enough to determine the sign of
		\[
		\phi_p(t)=(1+t)^p+(1-t)^p-2-(2^p-2)t^p.
		\]
		Since $\phi_p(0)=\phi_p(1)=0$, only $0<t<1$ needs consideration.
		For $0<t\le1$, put
		\[
		g_p(t)=\frac{(1+t)^p+(1-t)^p-2}{t^p}.
		\]
		For $0<t<1$, differentiation gives
		\begin{equation*}
			g_p'(t)=\frac{p\bigl(2-(1+t)^{p-1}-(1-t)^{p-1}\bigr)}{t^{p+1}}.
		\end{equation*}
		If $0<p<1$ or $p>2$, then $x^{p-1}$ is convex on $(0,\infty)$, so
		\[
		(1+t)^{p-1}+(1-t)^{p-1}\ge2.
		\]
		Thus $g_p'(t)\le0$, whence $g_p(t)\ge g_p(1)=2^p-2$ and $\phi_p(t)\ge0$.  If $1<p<2$, then $x^{p-1}$ is concave, and the same argument gives $g_p'(t)\ge0$, hence $g_p(t)\le g_p(1)=2^p-2$ and $\phi_p(t)\le0$.   For $p=1,2$, the expression is identically zero.
	\end{proof}
	
	We now integrate the scalar inequality in Lemma~\ref{lem:scalar} against the tracial joint spectral measure.  
	
	\begin{prop}\label{prop:hanner-two-point}
		Let $D,O\in \mathbb{H}_n$ satisfy
		$
		D\pm O\ge 0.
		$
		For $p>0$, define
		\[
		\Gamma_p(D,O):=\Tr(D+O)^p+\Tr(D-O)^p-2\Tr D^p-(2^p-2)\Tr |O|^p.
		\]
		Then
		\[
		\Gamma_p(D,O)
		\begin{cases}
			\ge0, & 0<p<1\ \hbox{or}\ p>2,\smallskip\\
			=0, & p=1,2,\smallskip\\
			\le0, & 1<p<2.
		\end{cases}
		\]
	\end{prop}
	
	\begin{proof}
		Since $D\pm O\ge0$, we have
		\[
		D=\frac{(D+O)+(D-O)}2\ge0.
		\]
		Let $\mu=\mu_{D,O}$ be the tracial joint spectral measure of $D,O$ determined in Lemma~\ref{thm:heinavaara}.  By Lemma~\ref{lem:support}, $a\ge |b|$ on the support of $\mu$, up to a null set.  Hence Lemma~\ref{lem:scalar} may be integrated over $\R^2$.
		
		Using \eqref{eq:tjsm-p}, and using $D\pm O\ge0$ and $D\ge0$, we have
		\begin{align*}
			\Tr(D\pm O)^p&=p(p+1)\int_{\R^2}|a\pm b|^p\,d\mu(a,b),\\
			\Tr D^p&=p(p+1)\int_{\R^2}a^p\,d\mu(a,b),\\
			\Tr|O|^p&=p(p+1)\int_{\R^2}|b|^p\,d\mu(a,b).
		\end{align*}
		Therefore
		\[
		\Gamma_p(D,O)=p(p+1)\int_{\R^2}\Phi_p(a,b)\,d\mu(a,b),
		\]
		where $\Phi_p(a,b)$ is defined by \eqref{eq:scalar-phi}.
		Thus, the sign follows immediately from Lemma~\ref{lem:scalar}.
	\end{proof}
	
	\subsection{The case \texorpdfstring{$p>0$}{p > 0} in Theorem~\ref{thm:main}}
	\begin{proof}[Proof of the case $p> 0$ in Theorem~\ref{thm:main}]
		For $p>0$ define
		\[
		\Delta_p(A,B)
		=
		\Tr(A+B)^p
		-
		\Tr A^p
		-
		\Tr B^p
		-
		(2^p-2)\Tr(A^{1/2}BA^{1/2})^{p/2}.
		\]
		Consider the $2n\times2n$ matrices
		$
		D=
		\begin{pmatrix}
			A&0\\
			0&B
		\end{pmatrix},
		O=
		\begin{pmatrix}
			0&A^{1/2}B^{1/2}\\
			B^{1/2}A^{1/2}&0
		\end{pmatrix}
		$ and 	$J=
		\begin{pmatrix}
			I&0\\
			0&-I
		\end{pmatrix}$.
		Then $D+O\ge 0$
		and $D-O=J(D+O)J\ge0$.  Hence Proposition~\ref{prop:hanner-two-point}
		applies to $D,O$.
		
		Since $D-O$ and $D+O$ are unitarily similar,
		\[
		\Tr(D-O)^p=\Tr(D+O)^p.
		\]
		Therefore Proposition~\ref{prop:hanner-two-point} gives the corresponding
		sign of
		\begin{equation}\label{eq:block-step-main}
			2\Tr(D+O)^p
			-
			2\Tr D^p
			-
			(2^p-2)\Tr|O|^p .
		\end{equation}
		
		We now identify the three trace terms in \eqref{eq:block-step-main}.  Let
		$
		W=
		\begin{pmatrix}
			A^{1/2}\\
			B^{1/2}
		\end{pmatrix}.
		$
		Then $D+O=WW^*$ and $W^*W=A+B$.  Hence $D+O$ and $A+B$ have the same
		non-zero eigenvalues, and since $p>0$,
		\begin{equation}\label{eq:block-main-1}
			\Tr(D+O)^p=\Tr(A+B)^p .
		\end{equation}
		Also,
		\begin{equation}\label{eq:block-main-2}
			\Tr D^p=\Tr A^p+\Tr B^p .
		\end{equation}
		Finally,
		\begin{equation}\label{eq:block-main-3}
			\Tr|O|^p
			=
			2\Tr(A^{1/2}BA^{1/2})^{p/2}.
		\end{equation}
		Substituting \eqref{eq:block-main-1}, \eqref{eq:block-main-2}, and
		\eqref{eq:block-main-3} into \eqref{eq:block-step-main} gives
		$
		2\Delta_p(A,B).
		$
		Consequently,
		\[
		\Delta_p(A,B)
		\begin{cases}
			\ge0, & 0<p<1\ \text{or}\ p>2,\\
			=0, & p=1,2,\\
			\le0, & 1<p<2.
		\end{cases}
		\]
		This proves the case $p>0$ of Theorem~\ref{thm:main}.
	\end{proof}
	
	\subsection{Log-majorization and matrix geometric mean}
	The following singular-value product inequality is standard in matrix analysis; see {\cite[p. 94, Theorem IV.2.5]{Bha97}}.
	\begin{lem}[{\cite{Bha97}}]
		\label{lem:singular-sum}
		Let \(X,Y\in \mathbb{M}_n\). Then for $p>0$,
		\begin{equation*}
			\sum_{j=1}^{n} s_j(XY)^p
			\leq
			\sum_{j=1}^{n} s_j(X)^p s_j(Y)^p,
		\end{equation*}
		where \(s_1(A)\geq\cdots\geq s_n(A)\geq0\) denote the singular
		values of \(A\).
	\end{lem}

	We begin with the following multiplicative estimate.
	
	\begin{lem}\label{lem:first-negative}
		Let $A,B\in \mathbb P_n^+$. Then for $p>0$,
		\begin{equation*}
			\Tr A^{-p}+\Tr B^{-p}\ge2\Tr(A^{1/2}BA^{1/2})^{-p/2}.
		\end{equation*}
	\end{lem}
	
	\begin{proof} Notice that
		\begin{align*}
			\Tr(A^{1/2}BA^{1/2})^{-p/2}&=\sum_{j=1}^n s_j(B^{-1/2}A^{-1/2})^p\\
			&\le \sum_{j=1}^n s_j(B^{-1/2})^p s_j(A^{-1/2})^p \qquad \text{ (By Lemma~\ref{lem:singular-sum})}\\
			&=  \sum_{j=1}^n s_j(B^{-p})^{1/2} s_j(A^{-p})^{1/2}\\
			&\le (\Tr B^{-p})^{1/2}(\Tr A^{-p})^{1/2}\qquad \text{ (By Cauchy-Schwarz inequality)}\\
			&\le \frac{\Tr A^{-p}+\Tr B^{-p}}2.
		\end{align*}
	\end{proof}
	
	The matrix geometric mean is defined by
	\begin{equation*}
		A\#B=A^{1/2}(A^{-1/2}BA^{-1/2})^{1/2}A^{1/2},
	\end{equation*}
	which is a central example in the Kubo-Ando theory of operator means \cite{KA80}; see also Ando, Li, and Mathias~\cite{ALM04} and Bhatia~\cite{Bha07}. 
	Let $A,B\in\mathbb{P}_n^+$, and arrange their eigenvalues in
	nonincreasing order:
	$
	\lambda_1(A)\geq\lambda_2(A)\geq\cdots\geq\lambda_n(A)>0,
	$
	and similarly for $B$.  We say that \emph{\(A\) is weakly log-majorized
		by \(B\)}, denoted by \(A\prec_{w\log}B\), if
	\[
	\prod_{j=1}^{k}\lambda_j(A)
	\leq
	\prod_{j=1}^{k}\lambda_j(B),
	\qquad k=1,\ldots,n.
	\]
	We say that \emph{\(A\) is log-majorized by \(B\)}, denoted by
	\(A\prec_{\log}B\), if, in addition,
	\(\prod_{j=1}^{n}\lambda_j(A)=\prod_{j=1}^{n}\lambda_j(B)\),
	or equivalently, \(\det A=\det B\).

	The celebrated Araki--Lieb--Thirring  log-majorization \cite{Ara90} states the following.
	\begin{lem}
		\label{lem:Araki-log-majorization}
		Let $A,B\in\mathbb{P}_n^{+}$. Then for $r\geq 1$,
		\begin{equation*}
			\bigl(A^{1/2}BA^{1/2}\bigr)^r
			\prec_{\log}
			A^{r/2}B^rA^{r/2}.
		\end{equation*}
	\end{lem}

	\begin{lem}\label{lem:second-negative}
		Let $A,B\in \mathbb P_n^+$. Then for $p>0$,
		\begin{equation*}
			\Tr(A+B)^{-p}\le 2^{-p}\Tr(A^{1/2}BA^{1/2})^{-p/2}.
		\end{equation*}
	\end{lem}
	
	\begin{proof}
		By matrix AM-GM inequality, we have $2A\# B\le A+B$.
		The eigenvalue monotonicity in the Loewner order gives
		\begin{equation*}
			\Tr(A+B)^{-p}\le 2^{-p}\Tr (A\# B)^{-p}.
		\end{equation*}
		It remains to prove
		\begin{equation}\label{eq:tr}
			\Tr (A\# B)^{-p}\le \Tr(A^{1/2}BA^{1/2})^{-p/2}.
		\end{equation}
		Let $
		C=\bigl(A^{-1/2}BA^{-1/2}\bigr)^{1/2}.
		$
		Then 
		$
		A\# B=A^{1/2}CA^{1/2}.
		$
		Applying Lemma~\ref{lem:Araki-log-majorization} to \(A\) and \(C\)
		with \(r=2\), we obtain
		\[
		(A\# B)^2
		=\bigl(A^{1/2}CA^{1/2}\bigr)^2
		\prec_{\log}
		AC^2A
		=A^{1/2}BA^{1/2}.
		\]
		Recall that if \(X\prec_{\log}Y\), then
		\(\operatorname{Tr}X^q\leq\operatorname{Tr}Y^q\) for every
		\(q\in\mathbb{R}\). 
		Taking \(q=-p/2\) in the above log-majorization yields \eqref{eq:tr}.
	\end{proof}
	
	Next, we prove the case $p<0$ in Theorem~\ref{thm:main}. 
	\begin{proof}[Proof of the case $p<0$ in Theorem~\ref{thm:main}]  Set $r=-p>0$.
		Notice that
		\[
		\begin{aligned}
			&\quad \Tr A^{-r}+\Tr B^{-r}+(2^{-r}-2)\Tr(A^{1/2}BA^{1/2})^{-r/2}\\
			&=\bigl(\Tr A^{-r}+\Tr B^{-r}-2\Tr(A^{1/2}BA^{1/2})^{-r/2}\bigr)+2^{-r}\Tr(A^{1/2}BA^{1/2})^{-r/2}\\
			&\ge2^{-r}\Tr(A^{1/2}BA^{1/2})^{-r/2} \qquad \text{ (By Lemma~\ref{lem:first-negative})}\\
			&\ge\Tr(A+B)^{-r}. \qquad\qquad\qquad \text{ (By Lemma~\ref{lem:second-negative})}
		\end{aligned}
		\]
	\end{proof}
	
	\begin{proof}[Proof of Theorem~\ref{thm:main}]
		
		We have already proved the cases \(p>0\) and \(p<0\). When $p=0$, we have
		\[
		\Tr(A+B)^0
		=
		\Tr I
		=
		\Tr A^0+\Tr B^0-(\Tr I),
		\]
		which is exactly equality in \eqref{eq:McC_refine}.
		Equality for $A=B>0$ follows by direct substitution for every real $p$.
		Increasing any one coefficient in a lower bound, or decreasing any
		one coefficient in an upper bound, makes that equality fail while
		the other two coefficients remain fixed. This proves the stated
		componentwise optimality and completes the proof.
	\end{proof}
	
	\section{Hanner's inequality below one}
	\label{sec:hanner-below-one}
	\begin{proof}[Proof of Proposition~\ref{prop:hanner-no-direction}]
		Define the Hanner defect by
		\[
		\Delta_p(A,B)
		=
		\|A+B\|_p^p+\|A-B\|_p^p
		-
		\Bigl[
		\bigl(\|A\|_p+\|B\|_p\bigr)^p
		+
		\bigl|\|A\|_p-\|B\|_p\bigr|^p
		\Bigr].
		\]
		
		First, take
		$
		A=
		\begin{pmatrix}
			1&0\\
			0&0
		\end{pmatrix},
		$ and $
		B=
		\begin{pmatrix}
			0&0\\
			0&1
		\end{pmatrix}.
		$
		Then
		\[
		A+B=I_2,
		\qquad
		A-B=\operatorname{diag}(1,-1),
		\]
		and hence
		\[
		\|A+B\|_p^p+\|A-B\|_p^p=2+2=4.
		\]
		On the other hand,
		$
		\|A\|_p=\|B\|_p=1,
		$
		so that
		\[
		\bigl(\|A\|_p+\|B\|_p\bigr)^p
		+
		\bigl|\|A\|_p-\|B\|_p\bigr|^p
		=
		2^p.
		\]
		Since $0<p<1$, we have $4>2^p$. Thus
		\[
		\Delta_p(A,B)>0.
		\]
		
		Conversely, let $P_n=\diag(1,0,\ldots,0)\in\M_n$ and take
		$
		A=P_n$ and 
		$	B=I_n.$
		Then
		\[
		A+B=\operatorname{diag}(2,1,\ldots,1),
		\qquad
		A-B=\operatorname{diag}(0,-1,\ldots,-1),
		\]
		and therefore
		\[
		\|A+B\|_p^p+\|A-B\|_p^p
		=
		2^p+2(n-1).
		\]
		Moreover,
		\[
		\|A\|_p=1,
		\qquad
		\|B\|_p=n^{1/p}.
		\]
		Put $t=n^{1/p}$. Then
		\[
		\begin{aligned}
			-\Delta_p(A,B)&= (1+n^{1/p})^p+(n^{1/p}-1)^p
			-
			\bigl(2^p+2(n-1)\bigr)  \\
			&=
			(t+1)^p+(t-1)^p-2t^p+(2-2^p).
		\end{aligned}
		\]
		Taylor's formula gives $(t+1)^p+(t-1)^p-2t^p=O(t^{p-2})$. Hence, as $t\to\infty$,
		\[
		(t+1)^p+(t-1)^p-2t^p\longrightarrow 0.
		\]
		Since $2-2^p>0$ for $0<p<1$, it follows that for all sufficiently large
		$n$,
		\[
		(1+n^{1/p})^p+(n^{1/p}-1)^p
		>
		2^p+2(n-1).
		\]
		Thus, for such $n$,
		\[
		\Delta_p(P_n,I_n)<0.
		\]
		
		Therefore for every $0<p<1$, the Hanner defect can be both positive and
		negative. Consequently, no fixed direction is possible.
	\end{proof}
	
	\begin{proof}[Proof of Theorem~\ref{thm:hanner-positive}]
		Set $\mathcal M=\bigl(\begin{smallmatrix}A&B\\B&A\end{smallmatrix}\bigr)$
		and $W=2^{-1/2}\bigl(\begin{smallmatrix}I_n&I_n\\I_n&-I_n\end{smallmatrix}\bigr)$.
		Then $W$ is unitary and
		$W^*\mathcal MW=\diag(A+B,A-B)\ge0$.
		Apply Theorem~\ref{thm:reverse-king-half-one}, proved in
		Section~\ref{sec:reverse-king}, to $\mathcal M$.
		Its scalar compression is
		$\bigl(\begin{smallmatrix}a&b\\b&a\end{smallmatrix}\bigr)\ge0$,
		where $a=\|A\|_p$ and $b=\|B\|_p$. Thus $a\ge b$, and the two
		eigenvalues of this compression are $a+b$ and $a-b$.
		Taking $p$th powers in \eqref{eq:reverse-king} gives
		\eqref{eq:hanner-positive}.
	\end{proof}
	
	\section{A reverse King inequality below the Banach range}
	\label{sec:reverse-king}

	If $q<0$ and $A\in \mathbb{P}_n^+$, we  use the standard negative-exponent notation
	\[
	\|A\|_q:=\bigl(\Tr A^q\bigr)^{1/q}.
	\]
	The purpose of this section is to prove Theorem~\ref{thm:reverse-king-half-one}.

	\subsection{Three auxiliary inequalities}
	
	We begin with a reverse H\"older inequality and the associated negative-exponent duality formula, which, in the range \(0<p<1\), serve as the natural analogues of the usual H\"older inequality and Schatten duality in the Banach setting; see \cite[Section~6, Proposition~6.1, Example~6.2, and Eq.~(6.8)]{BH14}; for the explicit optimizer, see also Chayes~\cite[Lemma~3.1]{Cha21b}. To make the exposition self-contained and to avoid any possible confusion, we include a complete proof.
	\begin{lem}
		\label{lem:negative-holder}
		Let $0<p<1$ and put
		$
		q=p/(p-1)<0.
		$
		Let  $A\in \mathbb P_n$ and $B\in \mathbb P_n^+$. Then
		\begin{equation}
			\Tr(AB)\geq \|A\|_p\,\|B\|_q.
			\label{eq:negative-holder}
		\end{equation}
		Consequently, for every $A\in \mathbb P_n^+$,
		\begin{equation}
			\|A\|_p
			=
			\inf_{\substack{B\in \mathbb P_n^+,\, \|B\|_q=1}}
			\Tr(AB),
			\label{eq:negative-duality}
		\end{equation}
		and the infimum is attained at
		$	B_A={A^{p-1}}/{\|A\|_p^{p-1}}.
		$
		In particular, the map $A\mapsto\|A\|_p$ is concave on the positive
		cone.
	\end{lem}

	\begin{proof}
		We first prove the scalar inequality. The zero case is immediate.
		By homogeneity, it is enough to
		consider nonnegative numbers $x_1,\dots,x_n$ and positive numbers
		$y_1,\dots,y_n$ satisfying
		$
		\sum_{i=1}^n x_i^p=1,
		$ and $
		\sum_{i=1}^n y_i^q=1.
		$
		Set
		$
		u_i=x_i^p,
		$ and $
		v_i=y_i^q
		$.
		Since $1/p+1/q=1$, we have, for every $i=1,\dots,n$,
		\[
		x_i y_i
		=
		v_i\left(\frac{u_i}{v_i}\right)^{1/p}.
		\]
		The function $t\mapsto t^{1/p}$ is convex. Jensen's inequality therefore
		gives
		\[
		\sum_{i=1}^n x_i y_i
		=
		\sum_{i=1}^n
		v_i\left(\frac{u_i}{v_i}\right)^{1/p}
		\geq
		\left(\sum_{i=1}^n u_i\right)^{1/p}
		=
		1.
		\]
		This is the scalar reverse H\"older inequality.
		
		Diagonalize
		\[
		A=U\diag(a_1,\dots,a_n)U^*,
		\qquad
		B=V\diag(b_1,\dots,b_n)V^*,
		\]
		and let $W=U^*V$. Then
		\[
		\begin{aligned}
			\Tr(AB)
			&=
			\Tr\!\left(
			\diag(a_1,\dots,a_n)\,
			W\diag(b_1,\dots,b_n)W^*
			\right) \\
			&=
			\sum_{i=1}^n\sum_{j=1}^n
			a_i b_j |W_{ij}|^2.
		\end{aligned}
		\]
		Since $W$ is unitary, the matrix
		$
		D:=\bigl(|W_{ij}|^2\bigr)_{i,j=1}^n
		$
		is doubly stochastic. By the Birkhoff theorem (see~\cite[Theorem~5.21]{Zha11}), $D$ is a
		convex combination of permutation matrices. Thus there exist
		permutations $\sigma_1,\dots,\sigma_N$ of $\{1,\ldots,n\}$ and coefficients
		$\lambda_1,\dots,\lambda_N\geq0$ satisfying
		$
		\sum_{k=1}^N\lambda_k=1
		$
		such that, for every $i,j=1,\dots,n$,
		\[
		|W_{ij}|^2
		=
		\sum_{k=1}^N
		\lambda_k\mathbf{1}_{\{j=\sigma_k(i)\}}.
		\]
		Consequently,
		\[
		\begin{aligned}
			\Tr(AB)
			&=
			\sum_{i=1}^n\sum_{j=1}^n
			a_i b_j
			\sum_{k=1}^N
			\lambda_k\mathbf{1}_{\{j=\sigma_k(i)\}} \\
			&=
			\sum_{k=1}^N\lambda_k
			\sum_{i=1}^n a_i b_{\sigma_k(i)}.
		\end{aligned}
		\]
		For every $k=1,\dots,N$, the scalar reverse H\"older inequality applied
		to the sequences $(a_i)_{i=1}^n$ and
		$(b_{\sigma_k(i)})_{i=1}^n$ gives
		\[
		\sum_{i=1}^n a_i b_{\sigma_k(i)}
		\geq
		\left(\sum_{i=1}^n a_i^p\right)^{1/p}
		\left(\sum_{i=1}^n b_{\sigma_k(i)}^q\right)^{1/q}.
		\]
		Since $\sigma_k$ is a permutation of $\{1,\dots,n\}$,
		\[
		\sum_{i=1}^n b_{\sigma_k(i)}^q
		=
		\sum_{i=1}^n b_i^q.
		\]
		Therefore,
		\[
		\begin{aligned}
			\Tr(AB)
			&\geq
			\sum_{k=1}^N\lambda_k
			\left(\sum_{i=1}^n a_i^p\right)^{1/p}
			\left(\sum_{i=1}^n b_i^q\right)^{1/q} \\
			&=
			\left(\sum_{k=1}^N\lambda_k\right)
			\left(\sum_{i=1}^n a_i^p\right)^{1/p}
			\left(\sum_{i=1}^n b_i^q\right)^{1/q} \\
			&=
			\|A\|_p\,\|B\|_q,
		\end{aligned}
		\]
		which proves \eqref{eq:negative-holder}.
		
		If $\|B\|_q=1$, then \eqref{eq:negative-holder} yields
		\[
		\Tr(AB)\geq\|A\|_p.
		\]
		The matrix $B_A$ satisfies
		\[
		\Tr(B_A^q)=1,
		\qquad
		\Tr(AB_A)=\|A\|_p,
		\]
		which proves \eqref{eq:negative-duality}. Finally,
		\eqref{eq:negative-duality} expresses the map
		$A\mapsto\|A\|_p$ as the infimum of a family of linear functionals,
		and hence this map is concave on the positive definite cone.
		Continuity extends concavity to the positive semidefinite cone.
	\end{proof}

	The next lemma is the only point at which the
	Araki--Lieb--Thirring inequality is used in this section.
	
	\begin{lem}
		\label{lem:gram-compression}
		Let $r\geq1$. Let $Z=[z_{ij}]_{i,j=1}^d\in \mathbb{P}_d^+$ and 
		$R_1,\dots,R_d\in \mathbb{M}_n$ such that
		$
		\|R_iR_i^*\|_r=1
		$
		for all $1\leq i\leq d$. Then
		\begin{equation}
			\left\|
			[z_{ij}R_iR_j^*]_{i,j=1}^d
			\right\|_r
			\leq
			\|Z\|_r.
			\label{eq:gram-compression}
		\end{equation}
	\end{lem}
	
	\begin{proof}
		Set
		\[
		\cD=\bigoplus_{i=1}^d R_i,
		\qquad
		B=Z\otimes I_n,
		\qquad
		Q=\cD B\cD^*.
		\]
		The nonzero spectra of $Q$ and
		$B^{1/2}\cD^*\cD B^{1/2}$ coincide. The
		Araki--Lieb--Thirring inequality (take traces in Lemma~\ref{lem:Araki-log-majorization}, extending to semidefinite matrices by continuity) gives
		\[
		\begin{aligned}
			\Tr Q^r &= \Tr (B^{1/2}\cD^*\cD B^{1/2})^r\\
			&\leq
			\Tr\!\left(B^r(\cD^*\cD)^r\right) \\
			&=
			\sum_{i=1}^d (Z^r)_{ii}
			\Tr(R_i^*R_i)^r \\
			&=
			\sum_{i=1}^d (Z^r)_{ii}
			=\Tr Z^r.
		\end{aligned}
		\]
		Taking $r$-th roots proves \eqref{eq:gram-compression}.
	\end{proof}
	
	We shall also need a Schatten estimate for the Kubo--Ando geometric mean.
	
	\begin{lem}
		\label{lem:geometric-mean}
		Let $P,Q\in \mathbb{P}_n^+$.
		If $G=P\#Q$, then
		\begin{equation}
			GP^{-1}G=Q,
			\qquad
			G^{-1}=P^{-1}\#Q^{-1}.
			\label{eq:geometric-identities}
		\end{equation}
		Moreover, for every $r>0$,
		\begin{equation*}
			\|P\#Q\|_r
			\leq
			\bigl(\|P\|_r\,\|Q\|_r\bigr)^{1/2}.
		\end{equation*}
	\end{lem}
	
	\begin{proof}
		The identities in \eqref{eq:geometric-identities} follow directly from
		the definition.
		The identity $GP^{-1}G=Q$ also gives $GQ^{-1}G=P$.
		Thus $U=P^{-1/2}GQ^{-1/2}$ is unitary, and
		$G=P^{1/2}UQ^{1/2}$.
		The multiplicative Schatten H\"older inequality, valid for all positive
		exponents, yields
		\[
		\|G\|_r
		\leq
		\|P^{1/2}\|_{2r}\,\|Q^{1/2}\|_{2r}
		=
		\bigl(\|P\|_r\,\|Q\|_r\bigr)^{1/2}.
		\]
	\end{proof}

	\subsection{Construction of a dual test matrix}
	
	We first consider a positive definite block matrix
	$M=\bigl(\begin{smallmatrix}X&Y\\Y^*&Z\end{smallmatrix}\bigr)$,
	where $X,Y,Z\in\M_n$. We may assume $Y\ge0$. Indeed, extend the
	partial isometry in $Y=V|Y|$ to a unitary $V\in\M_n$. Then
	\[
	\begin{pmatrix}V^*&0\\0&I_n\end{pmatrix}
	M\begin{pmatrix}V&0\\0&I_n\end{pmatrix}
	=\begin{pmatrix}V^*XV&|Y|\\|Y|&Z\end{pmatrix}>0.
	\]
	This conjugation preserves the Schatten functional of $M$ and of
	each of its blocks. We therefore work with $Y\ge0$ until the final
	approximation argument.
	
	The duality formula in Lemma~\ref{lem:negative-holder} gives explicit
	minimizers for $X$, $Z$, and their scalar compression $m_p$.
	We will use these three matrices to construct a dual test matrix for $M$.
	
	\begin{lem}\label{lem:normalized-data}
		Let $0<p<1$, let $q=p/(p-1)$, and suppose that
		$M=\bigl(\begin{smallmatrix}X&Y\\Y&Z\end{smallmatrix}\bigr)>0$
		with $Y\ge0$. Its scalar compression
		\[
		m_p=\begin{pmatrix}
			\|X\|_p&\|Y\|_p\\\|Y\|_p&\|Z\|_p
		\end{pmatrix}
		\]
		is positive definite. Define the scalar matrix $h=(h_{ij})_{i,j=1}^2$
		and the matrices $P,Q\in\M_n$ by
		\begin{equation}\label{eq:dual-data}
			h=\left(\frac{m_p}{\|m_p\|_p}\right)^{p-1},\qquad
			P=\left(\frac{X}{\|X\|_p}\right)^{p-1},\qquad
			Q=\left(\frac{Z}{\|Z\|_p}\right)^{p-1}.
		\end{equation}
		The matrix $h$ is real symmetric, $h_{12}\le0$, and
		\begin{equation}\label{eq:h-normalization}
			\|h\|_q=\|h^{-1}\|_{-q}=1,\qquad
			\Tr(m_ph)=\|m_p\|_p.
		\end{equation}
		The diagonal dual matrices satisfy
		\begin{equation}\label{eq:PQ-normalization}
			\begin{gathered}
				\|P\|_q=\|Q\|_q=\|P^{-1}\|_{-q}=\|Q^{-1}\|_{-q}=1,\\
				\Tr(XP)=\|X\|_p,\qquad \Tr(ZQ)=\|Z\|_p,
			\end{gathered}
		\end{equation}
		and their geometric mean satisfies
		\begin{equation}\label{eq:mean-normalization}
			\|P\#Q\|_q=\|(P\#Q)^{-1}\|_{-q}^{-1}\ge1.
		\end{equation}
	\end{lem}
	
	\begin{proof}
		\emph{Step 1: positivity of the scalar compression.}
		The Schur complement gives $Z-YX^{-1}Y>0$, so
		$\|X^{-1/2}YZ^{-1/2}\|_\infty<1$. The multiplicative Schatten
		H\"older inequality yields
		\[
		\begin{aligned}
			\|Y\|_p
			&\le \|X^{1/2}\|_{2p}\,
			\|X^{-1/2}YZ^{-1/2}\|_\infty\,
			\|Z^{1/2}\|_{2p}\\
			&<\bigl(\|X\|_p\|Z\|_p\bigr)^{1/2}.
		\end{aligned}
		\]
		Thus $m_p$ has positive diagonal entries and positive determinant.
		In particular, all the negative powers in \eqref{eq:dual-data} are
		well defined.
		
		\emph{Step 2: normalization and signs.}
		Since $q(p-1)=p$, we have
		$\Tr h^q=\|m_p\|_p^{-p}\Tr m_p^p=1$ and
		$\Tr(m_ph)=\|m_p\|_p^{1-p}\Tr m_p^p=\|m_p\|_p$.
		Also $(h^{-1})^{-q}=h^q$, which proves
		\eqref{eq:h-normalization}. Replacing $m_p$ by $X$ and by $Z$ in
		the same calculation proves \eqref{eq:PQ-normalization}.
		
		To check the sign of $h_{12}$, let $\lambda>\mu>0$ be the
		eigenvalues of $m_p$, when they are distinct. The spectral
		decomposition gives
		\[
		h_{12}=\|m_p\|_p^{1-p}\|Y\|_p
		\frac{\lambda^{p-1}-\mu^{p-1}}{\lambda-\mu}\le0.
		\]
		If the two eigenvalues coincide, $m_p$ is a scalar multiple of $I_2$
		and $h_{12}=0$. We will also use that every entry of $h^{-1}$ is
		nonnegative, as follows directly from
		\begin{equation}\label{eq:h-inverse-entries}
			h^{-1}=\frac1{h_{11}h_{22}-h_{12}^2}
			\begin{pmatrix}h_{22}&-h_{12}\\-h_{12}&h_{11}\end{pmatrix}.
		\end{equation}
		
		\emph{Step 3: the geometric mean.}
		Lemma~\ref{lem:geometric-mean} and \eqref{eq:PQ-normalization} give
		\[
		\|(P\#Q)^{-1}\|_{-q}
		=\|P^{-1}\#Q^{-1}\|_{-q}
		\le\bigl(\|P^{-1}\|_{-q}\|Q^{-1}\|_{-q}\bigr)^{1/2}=1.
		\]
		For every positive definite matrix, its $q$-functional is the
		reciprocal of the $(-q)$-functional of its inverse.
		Taking reciprocals proves \eqref{eq:mean-normalization}.
	\end{proof}
	
	We now replace the entries of $h$ by matrix blocks. The diagonal
	blocks use $P$ and $Q$, and the off-diagonal blocks use $P\#Q$.
	The sign $h_{12}\le0$ will turn the reverse H\"older lower bound
	for $\Tr(Y(P\#Q))$ into an upper bound for the full trace pairing.
	
	\begin{lem}\label{lem:lifted-certificate}
		For the block matrix $M$ in Lemma~\ref{lem:normalized-data}, let
		$h,P,Q$ be defined by \eqref{eq:dual-data}, and set
		\begin{equation}\label{eq:dual-block-matrix}
			H=\begin{pmatrix}
				h_{11}P&h_{12}(P\#Q)\\
				h_{12}(P\#Q)&h_{22}Q
			\end{pmatrix}.
		\end{equation}
		Then $H>0$, its inverse is
		\begin{equation}\label{eq:H-inverse}
			H^{-1}=\begin{pmatrix}
				(h^{-1})_{11}P^{-1}&(h^{-1})_{12}(P\#Q)^{-1}\\
				(h^{-1})_{12}(P\#Q)^{-1}&(h^{-1})_{22}Q^{-1}
			\end{pmatrix},
		\end{equation}
		and
		\begin{equation}\label{eq:MH-bound}
			\begin{aligned}
				0<\Tr(MH)
				&\le\|m_p\|_p+
				2h_{12}\|Y\|_p\bigl(\|P\#Q\|_q-1\bigr)\\
				&\le\|m_p\|_p.
			\end{aligned}
		\end{equation}
	\end{lem}
	
	\begin{proof}
		\emph{Step 1: positivity and inversion.}
		By \eqref{eq:geometric-identities},
		$(P\#Q)P^{-1}(P\#Q)=Q$. The Schur complement of $h_{11}P$
		in $H$ is therefore
		\[
		h_{22}Q-\frac{h_{12}^2}{h_{11}}(P\#Q)P^{-1}(P\#Q)
		=\left(h_{22}-\frac{h_{12}^2}{h_{11}}\right)Q>0,
		\]
		because $h>0$. Thus $H>0$. The same geometric-mean identity gives
		\[
		P^{-1}(P\#Q)Q^{-1}=(P\#Q)^{-1},\qquad
		(P\#Q)Q^{-1}(P\#Q)=P.
		\]
		Using these identities and the inverse formula
		\eqref{eq:h-inverse-entries}, direct block multiplication verifies
		that the right side of \eqref{eq:H-inverse} is the inverse of $H$.
		
		\emph{Step 2: the trace pairing.}
		By \eqref{eq:negative-holder},
		$\Tr(Y(P\#Q))\ge\|Y\|_p\|P\#Q\|_q$.
		Since $h_{12}\le0$, \eqref{eq:PQ-normalization} yields
		\[
		\begin{aligned}
			\Tr(MH)
			&=h_{11}\Tr(XP)+h_{22}\Tr(ZQ)+2h_{12}\Tr(Y(P\#Q))\\
			&\le h_{11}\|X\|_p+h_{22}\|Z\|_p
			+2h_{12}\|Y\|_p\|P\#Q\|_q\\
			&=\Tr(m_ph)+2h_{12}\|Y\|_p\bigl(\|P\#Q\|_q-1\bigr).
		\end{aligned}
		\]
		Now use \eqref{eq:h-normalization} and
		\eqref{eq:mean-normalization}. Finally,
		$\Tr(MH)=\Tr(M^{1/2}HM^{1/2})>0$, proving \eqref{eq:MH-bound}.
	\end{proof}
	
	\subsection{The range \texorpdfstring{$1/2\le p<1$}{1/2 <= p < 1}}
	
	When $p\ge1/2$, the positive exponent $-q=p/(1-p)$ is at least
	one. We can therefore apply Lemma~\ref{lem:gram-compression}
	to $H^{-1}$.
	
	\begin{prop}\label{prop:initial-interval}
		For $1/2\le p<1$ and every positive definite block matrix
		$M=\bigl(\begin{smallmatrix}X&Y\\Y^*&Z\end{smallmatrix}\bigr)$,
		the scalar compression $m_p$ in
		Theorem~\ref{thm:reverse-king-half-one} satisfies
		$\|M\|_p\le\|m_p\|_p$.
	\end{prop}
	
	\begin{proof}
		Make the unitary reduction preceding
		Lemma~\ref{lem:normalized-data}, and construct $h,P,Q$ and $H$
		as in \eqref{eq:dual-data} and \eqref{eq:dual-block-matrix}.
		In Lemma~\ref{lem:gram-compression}, take the scalar matrix
		$h^{-1}$, the exponent $-q$, and the two factors
	$
		P^{-1/2}$ and $(P\#Q)^{-1}P^{1/2}.
$
		Their products with their own adjoints are $P^{-1}$ and $Q^{-1}$,
		respectively, and the product of the first with the adjoint of the
		second is $(P\#Q)^{-1}$. Thus the block matrix in that lemma is
		exactly \eqref{eq:H-inverse}. Its diagonal factors have
		$(-q)$-functional one by \eqref{eq:PQ-normalization}. Hence
		\[
		\|H^{-1}\|_{-q}\le\|h^{-1}\|_{-q}=1,
		\qquad\text{and therefore}\qquad \|H\|_q\ge1.
		\]
		The reverse H\"older inequality \eqref{eq:negative-holder} and
		the trace bound \eqref{eq:MH-bound} now give
		\[
		\|M\|_p\le\frac{\Tr(MH)}{\|H\|_q}
		\le\Tr(MH)\le\|m_p\|_p.
		\]
	\end{proof}
	
	\subsection{A scalar estimate for the descent step}
	
	For smaller $p$, the exponent $p/(1-p)$ is below one.
	We will compare $H^{-1}$ with its scalar compression at this
	exponent. The compression differs from $h^{-1}$ only in that
	its off-diagonal entries are smaller. The following elementary
	estimate measures the resulting change in the Schatten functional.
	
	\begin{lem}\label{lem:scalar-deformation}
		Let $-1<q<0$ and let $K=(k_{ij})_{i,j=1}^2>0$ be real symmetric,
		with $k_{12}\ge0$ and $\Tr K^{-q}=1$. For $t\ge1$, define
		\[
		K_t=\begin{pmatrix}k_{11}&k_{12}/t\\k_{12}/t&k_{22}\end{pmatrix}.
		\]
		Then $K_t>0$ and
		\begin{equation}\label{eq:scalar-deformation}
			\bigl(1+2(t-1)(K^{1-q})_{12}(K^{-1})_{12}\bigr)
			\|K_t\|_{-q}\le1.
		\end{equation}
	\end{lem}
	
	\begin{proof}
		If $k_{12}=0$, then $K_t=K$, and equality holds in
		\eqref{eq:scalar-deformation}. Assume $k_{12}>0$.
		
		\emph{Step 1: concavity.}
		For $0\le x\le1$, define
		\[
		f(x)=\left\|\begin{pmatrix}
			k_{11}&xk_{12}\\xk_{12}&k_{22}
		\end{pmatrix}\right\|_{-q}.
		\]
		The matrix inside this functional is positive definite because its
		determinant is at least $\det K>0$. In particular $K_t>0$.
		Since $0<-q<1$, Lemma~\ref{lem:negative-holder} implies that $f$
		is concave. The normalization $\Tr K^{-q}=1$ and differentiation
		of the trace give $f(1)=1$ and $f'(1)=2k_{12}(K^{-q-1})_{12}$.
		The tangent-line inequality therefore gives
		\begin{equation}\label{eq:scalar-tangent}
			\|K_t\|_{-q}=f(1/t)\le1-(1-1/t)f'(1).
		\end{equation}
		
		\emph{Step 2: bounding the derivative.}
		Let $\lambda>\mu>0$ be the eigenvalues of $K$; they are distinct
		since $k_{12}>0$. For any real exponent, the off-diagonal entry
		of the corresponding power of $K$ is $k_{12}$ times the divided
		difference of that power at $\lambda$ and $\mu$. Also
		$(K^{-1})_{12}=-k_{12}/(\lambda\mu)$. Consequently,
		\[
		\begin{aligned}
			(K^{1-q})_{12}(K^{-1})_{12}
			&=-\frac{k_{12}^2}{\lambda\mu}
			\frac{\lambda^{1-q}-\mu^{1-q}}{\lambda-\mu}\\
			&\le k_{12}^2\frac{\lambda^{-q-1}-\mu^{-q-1}}{\lambda-\mu}
			=\frac{f'(1)}2\le0.
		\end{aligned}
		\]
		The middle inequality follows from
		\[
		\lambda^{1-q}-\mu^{1-q}
		+\lambda\mu(\lambda^{-q-1}-\mu^{-q-1})
		=(\lambda+\mu)(\lambda^{-q}-\mu^{-q})\ge0.
		\]
		Thus \eqref{eq:scalar-tangent} and $1-1/t\le t-1$ imply
		\[
		\|K_t\|_{-q}
		\le1-2(t-1)(K^{1-q})_{12}(K^{-1})_{12}.
		\]
		If the factor multiplying $\|K_t\|_{-q}$ in
		\eqref{eq:scalar-deformation} is nonpositive, the claimed bound
		is immediate. Otherwise multiply the last inequality by that
		factor. The resulting right side is
		$1-4(t-1)^2\bigl((K^{1-q})_{12}(K^{-1})_{12}\bigr)^2\le1$.
		This proves \eqref{eq:scalar-deformation}.
	\end{proof}
	
	\subsection{Descent in the exponent}
	
	We now prove the inequality at $p$ by applying the same inequality
	at the larger exponent $p/(1-p)$ to $H^{-1}$.
	
	\begin{prop}\label{prop:descent-step}
		Let $0<p<1/2$ and $q=p/(p-1)$. Suppose that
		Theorem~\ref{thm:reverse-king-half-one} holds at exponent $-q$ for
		every positive definite matrix with two equally sized blocks.
		Then it holds at exponent $p$ for every such matrix.
	\end{prop}
	
	\begin{proof}
		Reduce to $Y\ge0$ and construct $h,P,Q$ and $H$ from $p$ and $M$
		as in \eqref{eq:dual-data} and \eqref{eq:dual-block-matrix}.
		In this range $-1<q<0$. Put $t=\|P\#Q\|_q$, so $t\ge1$
		by \eqref{eq:mean-normalization}.
		
		\emph{Step 1: compress the inverse test matrix.}
		Every entry of $h^{-1}$ is nonnegative by
		\eqref{eq:h-inverse-entries}. Using \eqref{eq:H-inverse},
		\eqref{eq:PQ-normalization}, and
		$\|(P\#Q)^{-1}\|_{-q}=1/t$, the scalar $(-q)$-compression of
		$H^{-1}$ is
		\[
		K_t=\begin{pmatrix}
			(h^{-1})_{11}&(h^{-1})_{12}/t\\
			(h^{-1})_{12}/t&(h^{-1})_{22}
		\end{pmatrix}.
		\]
		This is precisely the matrix $K_t$ in
		Lemma~\ref{lem:scalar-deformation} with $K=h^{-1}$.
		The hypothesis at exponent $-q$ therefore gives
		\begin{equation}\label{eq:H-descent-bound}
			\frac1{\|H\|_q}=\|H^{-1}\|_{-q}\le\|K_t\|_{-q}.
		\end{equation}
		
		\emph{Step 2: apply the scalar estimate.}
		By \eqref{eq:h-normalization}, $\Tr(h^{-1})^{-q}=1$.
		Moreover, $(1-p)(1-q)=1$, so \eqref{eq:dual-data} gives
		\[
		(h^{-1})^{1-q}=\frac{m_p}{\|m_p\|_p},\qquad
		\bigl((h^{-1})^{1-q}\bigr)_{12}h_{12}
		=\frac{h_{12}\|Y\|_p}{\|m_p\|_p}.
		\]
		Substituting $K=h^{-1}$ in \eqref{eq:scalar-deformation} yields
		\begin{equation}\label{eq:descent-scalar-bound}
			\left(1+\frac{2h_{12}\|Y\|_p}{\|m_p\|_p}(t-1)\right)
			\|K_t\|_{-q}\le1.
		\end{equation}
		The factor in parentheses is positive: by \eqref{eq:MH-bound},
		its product with $\|m_p\|_p$ is at least $\Tr(MH)>0$.
		
		\emph{Step 3: finish the dual estimate.}
		The reverse H\"older inequality \eqref{eq:negative-holder}, the
		inverse bound \eqref{eq:H-descent-bound}, and the trace bound
		\eqref{eq:MH-bound} give
		\[
		\begin{aligned}
			\|M\|_p
			&\le\frac{\Tr(MH)}{\|H\|_q}
			\le\Tr(MH)\|K_t\|_{-q}\\
			&\le\bigl(\|m_p\|_p+2h_{12}\|Y\|_p(t-1)\bigr)\|K_t\|_{-q}
			\le\|m_p\|_p.
		\end{aligned}
		\]
		The last inequality is \eqref{eq:descent-scalar-bound} multiplied
		by $\|m_p\|_p$.
	\end{proof}
	
	\subsection{Completion of the proof}
	
	We now finish the proof of Theorem~\ref{thm:reverse-king-half-one}.
	Finitely many applications of Proposition~\ref{prop:descent-step}
	cover any given exponent in $(0,1)$.
	
	\begin{proof}
		\emph{Step 1: all positive definite block matrices.}
		Proposition~\ref{prop:initial-interval} proves the result for
		$p\in[1/2,1)$. Suppose it has been proved for every
		$p\in[1/N,1)$, where $N\ge2$.
		If $p\in[1/(N+1),1/N)$, then
		\[
		\frac{p}{1-p}\in
		\left[\frac1N,\frac1{N-1}\right)\subseteq[1/N,1).
		\]
		The induction hypothesis and Proposition~\ref{prop:descent-step}
		prove the result at $p$. Induction on $N$ covers $(0,1)$ because
		every $p>0$ belongs to $[1/N,1)$ for some finite $N$.
		Equivalently, each application reduces the reciprocal exponent by
		one: $1/(p/(1-p))=1/p-1$.
		
		\emph{Step 2: positive semidefinite block matrices.}
		Let $M\ge0$. For $\varepsilon>0$, apply Step 1 to
		$M+\varepsilon I_{2n}$, whose scalar compression is
		\[
		m_p(\varepsilon)=
		\begin{pmatrix}
			\|X+\varepsilon I_n\|_p&\|Y\|_p\\
			\|Y\|_p&\|Z+\varepsilon I_n\|_p
		\end{pmatrix}>0.
		\]
		Thus $\|M+\varepsilon I_{2n}\|_p\le\|m_p(\varepsilon)\|_p$.
		As $\varepsilon\downarrow0$, continuity of the eigenvalues and of
		$s\mapsto s^p$ on $[0,\infty)$ gives
		$m_p(\varepsilon)\to m_p\ge0$ and
		$\|M\|_p\le\|m_p\|_p$. This proves the theorem.
	\end{proof}
	
	At $p=1$, the scalar compression is positive semidefinite by the same
	factorization argument, and equality holds:
	$\|M\|_1=\Tr X+\Tr Z=\Tr m_1=\|m_1\|_1$.
	The multiplicative constant one in
	Theorem~\ref{thm:reverse-king-half-one} is sharp, since equality holds
	whenever $Y=0$.
	
	\section{Unitary orbits, unitarily invariant norms, and mixed norms}
	\label{sec:families}
	
	The proofs use the square-sum argument for the three results in
	Section~\ref{subsec:families}. The unitarily invariant norm and unitary
	orbit estimates use equality of the square sums. For complex
	interpolation we also need the corresponding inequality for a scalar
	contraction, because the adjoint of a rectangular isometry need not be
	an isometry.
	
	\begin{lem}\label{lem:isometry-square-sum}
		Let $V=(v_{ij})\in\M_{m,n}(\C)$ satisfy $V^*V\le I_n$ and let
		$A_1,\ldots,A_n\in\M_d$. If $B_i=\sum_{j=1}^n v_{ij}A_j$, then
		\begin{equation}\label{eq:contraction-square-sum}
			\sum_{i=1}^m|B_i|^2\le\sum_{j=1}^n|A_j|^2.
		\end{equation}
		If $V^*V=I_n$, then
		\begin{equation}\label{eq:isometry-square-sum}
			\sum_{i=1}^m|B_i|^2=\sum_{j=1}^n|A_j|^2.
		\end{equation}
	\end{lem}
	
	\begin{proof}
		Regard $\mathcal A=(A_1,\ldots,A_n)^T$ as a block column. Then
		$\mathcal B=(B_1,\ldots,B_m)^T=(V\otimes I_d)\mathcal A$, and
		\[
		\sum_{i=1}^m|B_i|^2
		=\mathcal A^*(V^*V\otimes I_d)\mathcal A
		\le\mathcal A^*\mathcal A
		=\sum_{j=1}^n|A_j|^2.
		\]
		Equality holds when $V^*V=I_n$.
	\end{proof}
	
	\subsection{Unitarily invariant norms}
	
	We use the power-function case of the Jensen and superadditivity
	inequalities of Hirzallah and Kittaneh.
	
	\begin{lem}[{\cite[Lemma~2.1]{HK08}}]
		\label{lem:Jensen_type_inequalities_uin}
		Let $T_1,\ldots,T_k\in\mathbb P_d$, and let $\unorm{\cdot}$ be
		any unitarily invariant norm on $\M_d$. For $r\ge1$,
		\begin{equation}\label{eq:uin-powers}
			\unorm{\sum_{j=1}^kT_j^r}
			\le\unorm{\left(\sum_{j=1}^kT_j\right)^r}
			\le k^{r-1}\unorm{\sum_{j=1}^kT_j^r}.
		\end{equation}
		Both inequalities reverse for $0<r\le1$.
	\end{lem}
	
	\begin{proof}[Proof of Theorem~\ref{thm:uin}]
		Put $S=\sum_i|B_i|^2=\sum_j|A_j|^2$.
		For $p\ge2$, applying \eqref{eq:uin-powers} with exponent $p/2$
		to each family gives
		\begin{align*}
			\unorm{\sum_i|B_i|^p}&\le\unorm{S^{p/2}}
			\le n^{p/2-1}\unorm{\sum_j|A_j|^p},\\
			\unorm{\sum_j|A_j|^p}&\le\unorm{S^{p/2}}
			\le m^{p/2-1}\unorm{\sum_i|B_i|^p}.
		\end{align*}
		These are the two assertions for $p\ge2$.
		For $0<p\le2$, Lemma~\ref{lem:Jensen_type_inequalities_uin}
		reverses both chains, giving the stated inequalities.
	\end{proof}
	
	The trace inequalities obtained from Theorem~\ref{thm:uin} extend
	to families in $S_p$, giving \cite[Theorem~1.8]{Zha26a} in full.
	To see this, choose finite-rank orthogonal projections $P_N$ such that
	$A_jP_N\to A_j$ and $B_iP_N\to B_i$ in $S_p$ for every $i,j$.
	Such projections exist by finite-rank approximation, since the
	families are finite. Equality of the square sums is preserved:
	\[
	\sum_i|B_iP_N|^2
	=P_N\left(\sum_i|B_i|^2\right)P_N
	=P_N\left(\sum_j|A_j|^2\right)P_N
	=\sum_j|A_jP_N|^2.
	\]
	For each $N$, the ranges of $P_N$, $A_jP_N$, and $B_iP_N$ span
	a common finite-dimensional subspace, and all these operators
	vanish on its orthogonal complement. Apply the finite-dimensional
	trace inequalities on this subspace and pass to the limit in $S_p$.
	
	\subsection{Unitary orbits}
	
	The same argument has a matrix-order version. The two inequalities
	below follow from the unitary orbit results of Bourin and Lee:
	the first by induction from their two-variable inequalities, and the
	second by taking equal weights in their Jensen inequality.
	
	\begin{lem}[{\cite[Corollaries~2.4, 3.2 and Theorem~3.1]{BL12}}]
		\label{lem:unitary_orbit_conv_conc}
		Let $T_1,\ldots,T_k\in\mathbb P_d$ and $r\ge1$.
		There are unitaries $W_1,\ldots,W_k,W\in\M_d$ such that
		\begin{equation}\label{eq:orbit-powers}
			\sum_{j=1}^kW_jT_j^rW_j^*
			\le\left(\sum_{j=1}^kT_j\right)^r
			\le k^{r-1}W\left(\sum_{j=1}^kT_j^r\right)W^*.
		\end{equation}
		For $0<r\le1$, there are unitaries for which both inequalities reverse.
	\end{lem}
	
	\begin{proof}[Proof of Theorem~\ref{thm:unitary-orbit}]
		Let $S=\sum_i|B_i|^2=\sum_j|A_j|^2$.
		For $p\ge2$, apply the first inequality in \eqref{eq:orbit-powers}
		to the $B$-family and the second to the $A$-family, both with
		exponent $p/2$. This gives unitaries $W_1,\ldots,W_m,W$ such that
		\[
		\sum_{i=1}^mW_i|B_i|^pW_i^*
		\le S^{p/2}\le n^{p/2-1}W\left(\sum_{j=1}^n|A_j|^p\right)W^*.
		\]
		Conjugating this chain by $W^*$ and setting $U_i=W^*W_i$
		proves \eqref{eq:orbit-forward}.
		Interchanging the two families proves \eqref{eq:orbit-backward},
		since the argument uses only equality of the square sums.
		For $0<p\le2$, use the reversed form of
		\eqref{eq:orbit-powers} and the same conjugations.
	\end{proof}
	
	\subsection{Complex interpolation and mixed norms}
	
	For a Banach space $E$, the notation $\ell_s^k(E)$ means $E^k$
	with norm $(\sum_j\|x_j\|_E^s)^{1/s}$ when $s<\infty$ and
	$\max_j\|x_j\|_E$ when $s=\infty$.
	For a compatible Banach couple $(X_0,X_1)$, write
	$[X_0,X_1]_\theta$ for its complex interpolation space,
	$0<\theta<1$. Compatibility means that both spaces embed continuously
	in a common Hausdorff topological vector space. We use the following
	standard interpolation facts, in the finite-sum form needed here.
	Here $\|T\|_{X\to Y}=\sup_{\|x\|_X\le1}\|Tx\|_Y$ denotes
	the operator norm of a bounded linear map $T:X\to Y$.
	
	\begin{lem}[{\cite[Theorems~4.1.2 and~5.6.3]{BL76}}]
		\label{lem:interpolation}
		Let $(X_0,X_1)$ and $(Y_0,Y_1)$ be compatible Banach couples.
		Suppose a linear map $T:X_0+X_1\to Y_0+Y_1$ satisfies
		$\|T\|_{X_j\to Y_j}\le M_j$ for $j=0,1$. Then
		\begin{equation}\label{eq:interpolation-bound}
			\|T\|_{[X_0,X_1]_\theta\to[Y_0,Y_1]_\theta}
			\le M_0^{1-\theta}M_1^\theta.
		\end{equation}
		For a Banach space $E$, $1\le s_0,s_1\le\infty$, and $k\in\N$,
		\begin{equation}\label{eq:finite-sum-interpolation}
			[\ell_{s_0}^k(E),\ell_{s_1}^k(E)]_\theta=\ell_s^k(E)
			\quad\text{isometrically},\qquad
			\frac1s=\frac{1-\theta}{s_0}+\frac\theta{s_1}.
		\end{equation}
	\end{lem}
	
	The Schatten exponent will remain fixed throughout the interpolation.
	The essential point is to establish the endpoints on the whole
	domain of a coefficient contraction. This makes them available for
	both $U$ and $U^*$.
	
	\begin{lem}\label{lem:contraction-endpoints}
		Let $1\le p<\infty$ and let $V\in\M_{m,n}(\C)$ satisfy
		$V^*V\le I_n$. Define
		$T_V(A_1,\ldots,A_n)=(\sum_jv_{ij}A_j)_{i=1}^m$.
		Then
		\begin{align}
			\|T_V\|_{\ell_1^n(S_p)\to\ell_1^m(S_p)}
			&\le m^{1/2},\label{eq:TU-l1}\\
			\|T_V\|_{\ell_\infty^n(S_p)\to\ell_\infty^m(S_p)}
			&\le n^{1/2},\label{eq:TU-linfty}\\
			\|T_V\|_{\ell_p^n(S_p)\to\ell_p^m(S_p)}
			&\le
			\begin{cases}
				m^{1/p-1/2},&1\le p\le2,\\
				n^{1/2-1/p},&2\le p<\infty.
			\end{cases}\label{eq:TU-diagonal}
		\end{align}
	\end{lem}
	
	\begin{proof}
		\emph{Step 1: the outer exponents $1$ and $\infty$.}
		Both the rows and the columns of $V$ have Euclidean norm at most
		one, since $V^*V\le I_n$ and $VV^*\le I_m$.
		If $B_i=\sum_jv_{ij}A_j$, the triangle inequality and
		Cauchy--Schwarz give
		\[
		\sum_i\|B_i\|_p
		\le\sum_j\left(\sum_i|v_{ij}|\right)\|A_j\|_p
		\le\sqrt m\sum_j\|A_j\|_p.
		\]
		Likewise,
		$\max_i\|B_i\|_p\le
		(\max_i\sum_j|v_{ij}|)\max_j\|A_j\|_p
		\le\sqrt n\max_j\|A_j\|_p$.
		These prove \eqref{eq:TU-l1} and \eqref{eq:TU-linfty}.
		
		\emph{Step 2: the outer exponent $p$.}
		First suppose the $A_j$ are matrices and put $r=p/2$,
		$S_A=\sum_j|A_j|^2$, and $S_B=\sum_i|B_i|^2$.
		By \eqref{eq:contraction-square-sum}, $S_B\le S_A$.
		Eigenvalue monotonicity gives $\Tr S_B^r\le\Tr S_A^r$ for every
		$r>0$. If $r\ge1$, the trace-norm instance of
		\eqref{eq:uin-powers} gives
		\[
		\sum_i\|B_i\|_p^p\le\Tr S_B^r
		\le\Tr S_A^r\le n^{r-1}\sum_j\|A_j\|_p^p.
		\]
		If $0<r\le1$, its reversed form gives
		\[
		\sum_i\|B_i\|_p^p\le m^{1-r}\Tr S_B^r
		\le m^{1-r}\Tr S_A^r
		\le m^{1-r}\sum_j\|A_j\|_p^p.
		\]
		Taking $p$th roots proves \eqref{eq:TU-diagonal} for matrices.
		For a finite family in $S_p$, choose common finite-rank orthogonal
		projections $E_N$ such that $E_NA_jE_N\to A_j$ in $S_p$ for every
		$j$. Apply the matrix estimate to these compressions and pass to
		the limit. The finite linear combinations defining $T_V$ also
		converge in $S_p$, which proves the general case.
	\end{proof}
	
	\begin{proof}[Proof of Theorem~\ref{thm:schatten_p_q_norms}]
		We prove the forward bound for every scalar contraction $V$, with
		$T_V$ as in Lemma~\ref{lem:contraction-endpoints}.
		The desired forward estimate is then the case $V=U$.
		All interpolation below uses \eqref{eq:interpolation-bound} and
		\eqref{eq:finite-sum-interpolation} with $E=S_p$.
		
		\emph{Step 1: $1\le p\le2$.}
		If $1\le s\le p$ and $p>1$, choose $\theta$ by
		$1/s=1-\theta+\theta/p$. Interpolating \eqref{eq:TU-l1} and
		the first case of \eqref{eq:TU-diagonal} yields
		\[
		\|T_V\|_{\ell_s^n(S_p)\to\ell_s^m(S_p)}
		\le m^{(1-\theta)/2+\theta(1/p-1/2)}
		=m^{1/s-1/2}.
		\]
		The case $p=s=1$ is already \eqref{eq:TU-l1}.
		If $p\le s\le\infty$, take $1/s=(1-\theta)/p$.
		Interpolation between \eqref{eq:TU-diagonal} and
		\eqref{eq:TU-linfty} gives
		\[
		\|T_V\|_{\ell_s^n(S_p)\to\ell_s^m(S_p)}
		\le m^{(1-\theta)(1/p-1/2)}n^{\theta/2}
		=m^{(1-p/2)/s}n^{(1-p/s)/2}.
		\]
		The endpoint choices $\theta=0,1$ refer directly to the corresponding
		endpoint bounds.
		
		\emph{Step 2: $2\le p<\infty$.}
		For $1\le s\le p$, use $1/s=1-\theta+\theta/p$ again, this time
		with the second case of \eqref{eq:TU-diagonal}. The result is
		\[
		\|T_V\|_{\ell_s^n(S_p)\to\ell_s^m(S_p)}
		\le m^{(1-\theta)/2}n^{\theta(1/2-1/p)}
		=m^{\frac{p/s-1}{2(p-1)}}
		n^{\frac{(p-2)(1-1/s)}{2(p-1)}}.
		\]
		For $p\le s\le\infty$, take $1/s=(1-\theta)/p$ and interpolate
		with \eqref{eq:TU-linfty}. This gives
		\[
		\|T_V\|_{\ell_s^n(S_p)\to\ell_s^m(S_p)}
		\le n^{(1-\theta)(1/2-1/p)+\theta/2}
		=n^{1/2-1/s}.
		\]
		Steps 1 and 2 prove
		$\|T_V\|_{\ell_s^n(S_p)\to\ell_s^m(S_p)}\le C_{p,s}(m,n)$,
		with the constants in \eqref{eq:mixed-constants}.
		
		\emph{Step 3: the reverse estimate.}
		Since $U^*U=I_n$, the matrix $U^*$ is a contraction and
		$T_{U^*}T_U=I$. Apply the bound just proved to $V=U^*$, with
		$m,n$ interchanged, and then to the tuple $B=T_UA$:
		\[
		\|A\|_{\ell_s^n(S_p)}
		=\|T_{U^*}B\|_{\ell_s^n(S_p)}
		\le C_{p,s}(n,m)\|B\|_{\ell_s^m(S_p)}.
		\]
		This is \eqref{eq:ell-s-backward}.
	\end{proof}

\end{document}